\documentclass[a4paper,11pt]{article}
\usepackage{amsmath,amscd,amsfonts,amssymb,epsf,latexsym}

\makeatletter

\long\def\@makefnt#1{\parindent 1em\noindent
             \hb@xt@1.8em{\hss\@textsuperscript{}}#1}
\long\def\@ftntext#1{\insert\footins{%
     \reset@font\footnotesize
     \interlinepenalty\interfootnotelinepenalty
     \splittopskip\footnotesep
     \splitmaxdepth \dp\strutbox \floatingpenalty \@MM
     \hsize\columnwidth \@parboxrestore
     \color@begingroup
       \@makefnt{%
         \rule\z@\footnotesep\ignorespaces#1\@finalstrut\strutbox}%
     \color@endgroup}}%
\def\subjclass#1{%
   \@ftntext{2000 {\itshape Mathematics Subject Classification.}\enspace
#1.}}
\def\keywords#1{%
   \@ftntext{{\itshape Key words and phrases.}\enspace #1.}}
\makeatother

\def\moins{\raise 1pt\hbox{{$\scriptstyle -$}}}
\def\plus{\raise 1pt\hbox{{$\scriptstyle +$}} }

\newtheorem{theorem}{Theorem}
\newtheorem{proposition}[theorem]{Proposition}
\newtheorem{lemma}[theorem]{Lemma}

\newtheorem{remark}[theorem]{Remark}
\newtheorem{definition}[theorem]{Definition}

\newtheorem{example}[theorem]{Example}

\def\qed{~\hbox{$\Box$}}

\def\Aut{\mathop{\rm Aut}}
\def\cJ{\mathop{\rm J}}

\def\codim{\mathop{\rm codim}}

\def\dim{\mathop{\rm dim}}
\def\Sym{\mathop{\rm Sym}}

\def\corank{\mathop{\rm corank}}

\def\cT{{\mathcal T}}
\def\cJ{{\mathcal J}}

\def\L{{\cal L}}
\def\C{{\mathbb C}}
\def\Z{{\mathbb Z}}
\def\Q{\widetilde{Q}}
\def\X{{\mathbb X}}

\def\s{\vskip6pt}

\def\ss{\vskip7pt}

\begin{document}

\title{\bf Positivity of Thom polynomials II: \\
the Lagrange singularities}

\author{Ma\l gorzata Mikosz\thanks{Research supported by the KBN grant NN201 387034.}\\
\small Warsaw University of Technology,\\
\small Pl.~Politechniki 1, 00-661, Warszawa, Poland\\
\small emmikosz@mech.pw.edu.pl \and
Piotr Pragacz\\
\small Institute of Mathematics of Polish Academy of Sciences\\
\small \'Sniadeckich 8, 00-956 Warszawa, Poland\\
\small P.Pragacz@impan.gov.pl \and
Andrzej Weber\thanks{Research supported by the KBN grant NN201 387034.}
\\
\small Department of Mathematics of Warsaw University\\
\small Banacha 2, 02-097 Warszawa, Poland\\
\small aweber@mimuw.edu.pl}



\date{}

\maketitle

\begin{abstract}
We study Thom polynomials associated with Lagrange singularities. We expand them
in the basis of $\Q$-functions. This basis plays a key role in the Schubert calculus of isotropic
Grassmannians. We prove that the $\Q$-function expansions of Thom polynomials of Lagrange singularities have
always nonnegative coefficients. This is an analog of a result on Thom polynomials
of mapping singularities and Schur S-functions, established formerly by the two last authors.
\end{abstract}

\section{Introduction}

In the present paper, we study the {\it Thom polynomials} of Lagrange
singularities.
This paper concerns the structure of these polynomials and is
a continuation of \cite{PW1}, where the case of singularities of maps
was investigated.
We look at the {\it positivity} properties of these Thom polynomials.
Such positivity properties are nowadays widely investigated in algebraic
geometry, see the monograph of Lazarsfeld \cite{L}.

\s

Let $L$ be a Lagrangian submanifold in the linear symplectic space
$$
V=W\oplus W^*\,,
$$
equipped with the standard symplectic form.

Classically, in real symplectic geometry, the {\it Maslov class}
(\cite{Ar1})
is represented by the cycle
$$
\Sigma=\{x\in L\;:\; \dim(T_xL\cap W^*)>0\}.
$$
This cycle is the locus of singularities of the projection $L\to W$.
Its cohomology class is integral and modulo 2 it is equal to $w_1(T^*L)$,
the
first Stiefel-Whitney class. In complex symplectic geometry the same
construction applied for a Lagrangian submanifold $L$ contained in a
symplectic manifold fibering over a base $B$ with Lagrangian fibers,
leads to the cohomology class which is equal to
$$
c_1(T^*L-T^*B)\,.
$$
The gene\-ra\-lizations of the Maslov class are Thom polynomials
associated with the higher order types of singularities. These
types are defined by imposing conditions on the higher order jets
of $L$, see Definition \ref{types}.

For real singularities, the associated cohomology classes were
studied by Arnold and Fuks (see, e.g., \cite{Fk}), Vassiliev
\cite{Va1}, Audin \cite{Au2}, and others. The complex case was
studied by Kazarian \cite{Ka2} who also computed a substantial
number of examples. The Thom polynomials in this case can
always be written as polynomials in the Chern classes of $T^*L \to T^*B$
(see Remark \ref{rel}).

\s

Given a Thom polynomial, one can expand it in different bases. The Thom
polynomial of a singularity of map, in general, is not a positive
combination
of monomials in Chern classes. As shown in \cite{PW1}, such a Thom
polynomial
is always a positive $\Z$-linear combination of Schur $S$-functions.

The Thom polynomial of a Lagrange singularity is, in general, neither a
positive
combination of monomials in Chern classes nor a positive combination of
Schur
$S$-functions.

In the present paper, we use $\Q$-functions of \cite{PR} and show
that the Thom polynomial of a Lagrange singularity, expanded in the basis
of these
$\Q$-functions, has nonnegative coefficients.

\s

Here is a brief outline of the content of the paper.

In Section \ref{lag}, we recall definitions and properties of Lagrange
singularities. We introduce the space of jets of Lagrange submanifolds
and define the notion of a ``Lagrange singularity type''.

In Section \ref{Q-LG}, we recall the algebraic properties of
$\Q$-functions
of \cite{PR}. We recall also their cohomological interpretation in terms
of Schubert classes for the Lagrangian Grassmannian  from \cite{P1}.

In Section \ref{thom}, we attach the ``Thom polynomial'' to a Lagrange
singularity
type. Then we state and prove our main result, Theorem \ref{main},
asserting that
the Thom polynomial of a Lagrange singularity type is a nonnegative
$\Z$-linear combination
of $\Q$-functions.

The proof of the theorem is quite different from the one in \cite{PW1}.
There are two reasons for that. First, the $\Q$-functions do not admit --
up to now -- a
characterization similar to that of Schur $S$-functions from \cite{FL1}.
Second, the Lagrangian
case is ``more rigid'' than the one for singularities of maps, and does
not admit a functorial
interpretation like that in \cite{PW1}.
The proof in the present paper relies on the computation of some normal
bundle in the space
of jets of Lagrangian submanifolds, see Lemma \ref{normal}. This lemma is
based on a result
about actions of linear transformations on jets of functions (see
Proposition \ref{tang})
which seems to be of independent interest. In the key step, the proof uses
deformation to the normal cone.

We finish Section \ref{thom} with some discussion of Thom polynomials
of Legendre singularities.

\s

In the last section, we list the $\Q$-function expansions of Thom polynomials
up to codimension 6.

\section{Jets of Lagrangian submanifolds}\label{lag}

Let us fix a positive integer $n$.
Suppose that $W$ be a complex linear space, where $\dim W=n$. Let
$$
V=W\oplus W^*
$$
be a linear symplectic space, equipped with the standard
symplectic form $\langle -,- \rangle$, defined by
$$
\langle (w_1,f_1),(w_2,f_2)\rangle=f_1(w_2)-f_2(w_1)
$$
for $w_i \in W$ and $f_i \in W^*$, $i=1,2$. We shall view $V$ as a symplectic
manifold. Writing $q=(q_1,\ldots,q_n)$ for the coordinates of $W$ and $p=(p_1,\ldots,p_n)$
for the dual coordinates of $W^*$, the symplectic form on $V$ is \ $\sum_{i=1}^n dp_i \wedge dq_i$.

\s

Denote by \ $\rho: V \to W$ the projection.

\s

Any germ $L$ of a Lagrangian submanifold of $V$ with the nonsingular
projection
$\rho_{|L}$ is a graph of a 1-form \ $\alpha: W\to W^*$.
The condition that $L$ is Lagrangian is equivalent to
  \ $d\alpha=0$.
Since we deal with germs (in fact, with their jets), we can write
  \ $\alpha=dF$ \
for some function $F:W \to {\C}$.

\s

In the present paper, we shall investigate (the germs of) {\it singular}
Lagrangian
submanifolds $L$, that is $L$ itself is smooth but the projection \
$\rho_{|L}$
is singular.

Here is the simplest example. Suppose that $\dim W=1$, and set
\begin{equation}\label{Adwa}
L=\{(q,p)\in W\oplus W^*\,:\,q=p^2 \}\,.
\end{equation}
(In the classification of Lagrange singularities, this corresponds to the
singularity of type $A_2$.)

To classify the (germs of) Lagrangian submanifolds, one introduces a
suitable notion
of a {\it generating family of a Lagrangian germ} (for motivation,
see \cite[Example 5 in \S18.3]{AGV}, and for a precise definition, see
\cite[\S19.1]{AGV}).
(These generating families are usually versal deformations of suitable
functions
$f: \C^{\bullet} \to \C$  \ ({\it loc.cit.})).
Then one introduces, in \cite[\S19.4]{AGV}, the notions of $R$-equivalence
and (stable)
$R^+$-equivalence of two generating families of Lagrangian germs.
The crucial result (see \cite[Theorem in \S19.4]{AGV} p.~304) asserts that
the germs of two
Lagrangian submanifolds are {\it Lagrangian equivalent} if and only if the
corresponding
{\it generating families} are stably $R^+$-equivalent. We illustrate these
issues by the
following example.

\begin{example} \rm The Lagrangian submanifold corresponding to the
singularity $A_3$,
i.e.~$f(x)=x^4$, is obtained in the following way. The generating family
$$
F(x,q)=x^4+q_1x^2+q_2x
$$
is a {\it universal deformation} of $f(x)=F(x,0)=x^4$ (with additional
condition $F(0,q)=0$).
The generating family defines the Lagrangian submanifold $L$ in the
following way
$$
L=\left\{(q,p)\in W\oplus W^*\,:\,\exists \ x \in\C\,,\;
\frac{\partial F(x,q)}{\partial x}=0\,,\;\frac{\partial
F(x,q)}{\partial q}=p\,\right\}\,,
$$
that is,
$$
4x^3+2q_1x+q_2=0\,, \ \ \ \ \ x^2=p_1\,, \ \ \ \ \ x=p_2\,.
$$
We thus obtain
$$
L=\left\{(q,p)\in W\oplus W^*\,:\,p_1=p_2^2\,,\;q_2=-(4p_2^3+2q_1p_2)\right\} \,.
$$
\end{example}

\s

We will study the space of germs of Lagrangian submanifolds
 $L\subset V$ passing through 0. This space has infinite
 dimension, which is inconvenient from the point of view of
 algebraic geometry. Therefore we fix once for all a nonnegative integer $k$
 and we identify two germs if they have the degree of tangency at
 0 bigger than $k$. This way we obtain the space of $k$-jets of
 Lagrangian manifolds denoted by $\L(V)$. This space is homogeneous.
 Every germ of a Lagrangian submanifold can be obtained from the
``distinguished'' Lagrangian submanifold $W$ by application of a
germ symplectomorphism preserving 0. We have the following
description
\begin{equation}
\L(V)=\Aut(V)/P\,,
\end{equation}
where $\Aut(V)$ is the {\it group of $k$-jet symplectomorphisms
preserving $0\in V$}, and $P$ is the stabilizer of $W$.

Denote now by \ $LG(V)$ \ the {\it Lagrangian Grassmannian}
parame\-trizing
all linear Lagrangian subspaces of $V$. This manifold is embedded in
$\L(V)$ in a natural way.

On the other hand, we have the {\it Gauss map} \ $\pi: \L(V)\to LG(V)$,
which is a retraction to $LG(V)$, defined for a Lagrangian submanifold $L$
by
$$
\pi(L)=T_0(L)\,,
$$
the tangent space of $L$ at $0\in L$.

\begin{lemma}\label{fib} The fiber of the projection $\pi$ is the affine
space isomorphic
to \ $\bigoplus_{i=3}^{k+1}{{\Sym}^i}(W^*)$.\footnote{We point out that \ $\pi: \L(V)\to LG(V)$ \
is {\it not} a vector bundle. This projection is of the form
$$
\L(V)=Sp(V)\times_{P'}\L(V)_W\to Sp(V)/P'=LG(V)\,,
$$
where $P'=P\cap Sp(V)$ is the group of linear symplectomorphisms
stabilizing $W$. Calculation shows that $P'$ {\it does not} act linearly
on $\L(V)_W$ already for $k=3$.}
\end{lemma}
{\it Proof.} The fiber $\L(V)_{W}=\pi^{-1}(W)$ consists of those (jets of)
Lagrangian submanifolds which have the tangent space equal to $W$.
Every Lagrangian submanifold $L$ with nondegenerate projection
onto $W$ is the graph of the differential of a function \ $F: W\to \C$ \
(note that $dF$ acts from $W$ to $W^*$). The condition: \ $0\in L$ \
corresponds to the condition: \ $dF(0)=0$, and the condition:
  \ $T_{0}(L)=W$ \ corresponds to vanishing of the second derivatives of
$F$ at $0$.
This proves the lemma.
\qed

\s

We end this section with a definition of a {\it Lagrange
singularity type}. Consider the following subgroup $H$ of
$\Aut(V)$ consisting of holomorphic symplectomorphisms preserving
the fibration $\rho: V\to W$.
Such symplectomorphisms are compositions of maps induced by biholomorphisms
of the base and differentials of functions on the base (see
\cite[\S18.5, Theorem, p.~284]{AGV}).

We say, following \cite[\S18.6]{AGV}, that two jets of Lagrangian
submanifolds are {\it Lagrangian equivalent}, if they belong to
the same orbit of $H$.

We study not only  individual orbits of the group $H$ but also
families of orbits. Also, together with an orbit we are forced to
consider its closure. Therefore we introduce the following definition.
\begin{definition}\label{types}
A {\it Lagrange singularity type} $\Sigma$ is any closed pure-dimensional
algebraic subset of $\L(V)$ which is invariant with respect to the action
of $H$.
\end{definition}
(In other words, a Lagrange singularity type is a closed algebraic
set which is the union of some Lagrangian equivalence classes.)

For instance, the closure of the orbit of the singularity
$A_2$ is described by the condition: $\corank (D\rho_{|L}(0)) \ge 1$.
Similarly, the closure of the orbit of the singularity $D_4$ is given
by the condition: $\corank (D\rho_{|L}(0)) \ge 2$.
The singularity class $P_8$ is not the closure of a single orbit. The
family $P_8$ has one parameter (i.e.~it has modality 1). An orbit
of a germ belonging to $P_8$ has codimension 7, while the
singularity type $P_8$ has codimension 6. It can be described by the condition:
$\corank (D\rho_{|L}(0)) \ge 3$.
To define the singularity type $A_3$ one has to consider
the degeneracy locus $S\subset L$ of the differential
$D\rho_{|L}$. The singularity type $A_3$ consists of the jets of $L$ for
which $S$ is singular or $\rho_{|S}$ is not an immersion. These
conditions can be translated into algebraic equations in $\L(V)$.

In the literature (\cite[\S21.3]{AGV}), one can find a notion of {\it stable} Lagrangian
equivalence classes. They are classified in small codimensions ({\it loc.cit.}).

\section{Lagrangian Grassmannians and $\Q$-functions}\label{Q-LG}

We start with recollections on $\Q$-functions of \cite{PR}.

\s

Let $\X$ be an alphabet\footnote{By an {\it alphabet} we understand a finite multiset of elements
in a commutative ring.}.
By $\X^2$ we shall denote the alphabet consisting of squares of elements
of $\X$.
Given an alphabet of variables $\X$, we shall
denote by $\Sym(\X)$ the ring of symmetric functions in $\X$. Given any alphabet $\X=\{x_1,x_2,\ldots\}$, we set
\begin{equation}
\Q_i(\X)=e_i(\X)=\sum_{j_1<\cdots <j_i} x_{j_1}\cdots x_{j_i}\,,
\end{equation}
the $i$th elementary symmetric function in $\X$\,.

Given two nonnegative integers $i \ge j$, we put
\begin{equation}
\Q_{i,j}(\X)=\Q_i(\X) \Q_j(\X)
+2\sum\limits^j_{p=1}(-1)^p\Q_{i+p}(\X) \Q_{j-p}(\X)\,.
\end{equation}

For example, we have \ $\Q_{i,i}(\X)= e_i(\X^2)$.

Given any partition \ $I=(i_1 \ge \cdots \ge i_h \ge 0)$,
where we can assume $h$ to be even, we set
\begin{equation}
\Q_I(\X) = \hbox{Pfaffian} (M)\,,
\end{equation}
where $M=(m_{p,q})$ is the $h\times h$ skew-symmetric matrix with
$$
m_{p,q}= \Q_{i_p,i_q}(\X)
$$
for $1 \le p < q \le h$.

For an alphabet $\X=\{x_1,x_2,\ldots \}$ of degree 1 variables, the degree
of $\Q_I(\X)$ is equal to $|I|:=i_1+\cdots+i_h$.

For a fixed positive integer $n$, let $X_n$ be an alphabet of $n$
variables of degree $1$. Then the set $\bigl\{\widetilde
Q_I(\X_n)\bigr\}$ indexed by all partitions such that $i_1\le n$
forms an additive basis of the ring $\Sym(\X_n)$ (see \cite{PR}).
We shall say that a partition is {\it strict} if its parts are
distinct. Then the set $\bigl\{\widetilde Q_I(\X_n)\bigr\}$
indexed by all strict partitions such that $i_1\le n$ forms a
basis of the ring $\Sym(\X_n)$ as a free $\Sym(\X_n^2)$-module
({\it loc.cit.}). The same assertions hold for a countable
alphabet of variables without restriction on $i_1$.

Let $c_1, c_2,\ldots$ be a sequence of commuting variables, where
$\deg(c_i)=i$,
and let $\X=\{x_1,x_2,\ldots\}$ be an alphabet of degree $1$ variables.
We get a ring isomorphism
$$
\Phi: \Sym(\X) \to {\Z}[c_1,c_2,\ldots]\,,
$$
by setting \ $\Phi(e_i(\X))=c_i$ for $i=1,2,\ldots$.

Given a partition $I$, we put
\begin{equation}\label{QI}
\Q_I=\Phi (\Q_I(\X))\,.
\end{equation}

If $E$ is a vector bundle, then we define
\begin{equation}
\Q_I(E)=\Q_I(\X)\,,
\end{equation}
where $\X$ is the alphabet of the {\it Chern roots} of $E$.
In other words, $\Q_I(E)$ is equal to $\Q_I$, where $c_i$
is specialized to $c_i(E)$, the $i$th Chern class of $E$,
$i=1,2,\ldots$.

\begin{remark}\label{link}\rm
The family of $\Q$-functions was invented and investigated in \cite{PR}
on the occasion of study of Lagrangian degeneracy loci. It is modelled on
the classical {\it Schur $Q$-functions} (see, e.g., \cite{P1}). More
precisely,
for a strict partition $I$, the Schur $Q$-function of a vector bundle $E$
is
obtained from $\Q_I$ by the substitution
$$
Q_I(E)=\Q_I(E-E^*)\,.
$$
See also Remark 5.3 in \cite{LP} for another link between these two
families
of functions. We refer the reader to \cite{PR} and \cite{LP} for detailed
studies
of $\Q$-functions.
\end{remark}

We shall now use $\Q$-functions to describe some cohomological
properties of the Lagrangian Grassmannian $LG(V)$.
First, we recall presentation of the cohomology ring of $LG(V)$
by generators and relations, that goes back to Borel \cite{Bo}.

\begin{proposition}\label{gr}
With the above notation, we have
\begin{equation}\label{koh}
H^*(LG(V), {\Z}) \cong
{\Z}[c_1,c_2,\ldots, c_n]/(\Q_{i,i})_{i=1,2,\ldots, n}\,.
\end{equation}
\end{proposition}
(Here the $c_i$'s correspond to the Chern classes of the dual of
the tautological subbundle on $LG(V)$.)

The Lagrangian Grassmannian $LG(V)$ has an algebraic cell decomposition
which is a particular case of classical Schubert-Bruhat cell
decomposition,
but here admits the following concrete form. Suppose that a general flag
$$
V_{\bullet}: \ V_1 \subset V_2 \subset \cdots \subset V_n \subset V
$$
of isotropic subspaces with $\dim \ V_i = i$, is given
(i.e., equivalently, $V_n$ is Lagrangian).
Given a partition $I=(n \ge i_1 > \cdots > i_h > 0)$, we define
\begin{equation}
\Omega_I(V_{\bullet})=
\{L\in LG(V): \ \dim\bigl(L\cap V_{n+1-i_p}\bigr)\ge p,
\ p=1,\ldots,h \}.
\end{equation}
-- a {\it Schubert variety} associated with $I$. Note that we have
\begin{equation}
\codim\bigl(\Omega_I(V_{\bullet}),LG(V)\bigr)=|I|
\,,
\end{equation}
and that the cohomology class of $\Omega_I(V_{\bullet})$, denoted by
$\Omega_I$,
does not depend on the choice of the flag $V_{\bullet}$.
Let us recall the following expression for $\Omega_I$ in terms of
$\Q$-functions.

\begin{theorem}\label{G} (\cite[Sect.6]{P1})
Let $V$ be a $2n$-dimensional linear symplectic space.
Then in $H^{2|I|}(LG(V), \Z)$, we have
\begin{equation}
\Omega_I =\Q_I(R^{*})\,,
\end{equation}
where $R$ is the tautological subbundle on $LG(V)$.
\end{theorem}
(In \cite[Sect.6]{P1}, this result was given in terms of the {\it special
Schubert classes} $\Omega_i$ but $\Omega_i=c_i(R^*)$.)

\s

Since the Schubert varieties are closures of the cells of a cellular decomposition
of $LG(V)$, the Schubert classes \{$\Omega_I$\}, $I$ strict and $i_1\le n$,
form a $\Z$-basis of $H^*(LG(V), \Z)$. Hence also the polynomials $\{\Q_I(R^*)\}$
indexed by the partitions from the same set, have this property.

\s

In the proof of Theorem \ref{main}, we shall use the following result.
For a strict partition $I=(n \ge i_1 > \cdots > i_h > 0)$, we denote by
$I'$
the strict partition whose set of parts complements $\{i_1,\ldots,i_h\}$
in
$\{1,\ldots,n\}$.

\begin{proposition}\label{dual}
For a strict partition $I$ with $i_1\le n$, there exists only one strict
partition $J$ with $j_1\le n$ and $|J|=\dim LG(V)-|I|$ , for which
$$
\Q_I(R^*)\cdot \Omega_J \ne 0\,.
$$
In fact, this $J$ is equal to $I'$, and we have
\begin{equation}
\int_{LG(V)} \Q_I(R^*)\cdot \Omega_{I'}=1.
\end{equation}
\end{proposition}
By virtue of Theorem \ref{G}, this proposition follows, e.g., from \cite{PR},
Theorem 5.23.
See also Example 4.2 (5) in \cite{P1}.

\section{Thom polynomials of Lagrange singularities and
$\Q$-functions}\label{thom}

A Lagrange singularity type \ $\Sigma\subset\L(V)$ defines the cohomology
class
$$
[\Sigma]\in H^*(\L(V), {\Z})\cong H^*(LG(V), {\Z})\,.
$$
Suppose that this class is equal to
$$
\sum_I \alpha_I \Q_I(R^*)\,,
$$
where the sum runs over strict partitions $I$ with $i_1\le n$\,, and
$\alpha_I \in \Z$ (it is important here to use the bundle $R^*$). Then
\begin{equation}\label{alpha}
\cT^\Sigma :=\sum_I \alpha_I \Q_I
\end{equation}
is called the {\it Thom polynomial} associated with the Lagrange
singularity type $\Sigma$.

\begin{example}\label{1} \rm
We list here the $\Q$-functions expansions of Thom polynomials
of some Lagrange singularities. They were computed in \cite{Ka2}
in the basis of monomials in Chern classes.

\bigskip

{\parindent 0pt

\s

$A_2$: \ ${\Q_1}$

\s

$A_3$: \ ${3\Q_2}$

\s

$D_4$: \ $\Q_{21}$

\s

$D_6$: \ ${12 \Q_{32}+  24\Q_{41}}$

\s

$A_7$: \ ${135 \Q_{321}+ 1275 \Q_{42}+ 2004\Q_{51}+  2520 \Q_{6}}$

\s

$P_8$: \ $\Q_{321}$.}

\s

Note that that all the coefficients in the formulas of Example
\ref{1} are nonnegative. For more extensive list of examples, see Section \ref{examples}.

The Thom polynomials of the singularities $A_2$,
$D_4$ and $P_8$ are equal to single $\Q$-functions because they
are defined by conditions involving only the differential of $\rho_{|L}$,
that is, by conditions defining single Schubert varieties in the Lagrangian
Grassmannian.

For example, $P_8$ is defined by the condition: $\corank(D\rho_{|L}(0))\ge 3$.
This singularity type is the closure of one dimensional family of
orbits. Therefore the germs of the type $P_8$ are {\it not} stable in
the sense of \cite[\S21.3]{AGV}. On the other hand, $P_8$ is worth mentioning
since it is the first example of modality that appears in the
classification. The Thom polynomial $\Q_{321}$ is an obstruction to
avoid this singularity by deformation.

\end{example}

\s

We state now our main result.

\begin{theorem}\label{main}
For any Lagrange singularity type $\Sigma$, all the coefficients
$\alpha_I$ in (\ref{alpha}) are nonnegative.
\end{theorem}

For the proof of the theorem, we need several preliminary results.

\s

First of all, we shall use the {\it nonnegativity property}
of globally generated bundles \cite{FL0} (see also \cite{F}, \cite{L}).
Let $E$ be a vector bundle on a variety $X$. By a {\it cone} in $E$,
we mean a subvariety of $E$ which is stable under the natural
${\mathbb G}_m$-action on $E$.
If $C\subset E$ is a cone, then one may intersect its cycle $[C]$ with
the zero-section of the vector bundle:
\begin{equation}
z(C,E):=s_E^*([C])\,,
\end{equation}
where \ $s_E^*: H^*(E, \Z)\to H^*(X, \Z)$ is the map
of the cohomology groups determined by the zero-section
$s_E: X \to E$.

We now record the following variant of a result from \cite{FL0}.

\begin{lemma}\label{int}
Let $\pi: E\to X$ be a globally generated bundle on
a proper homogeneous variety $X$. Let $C$ be a cone in $E$,
and let $Z$ be any algebraic cycle in $X$ of the complementary
dimension. Then the intersection \ $[C]\cdot [Z]$ is nonnegative.
\end{lemma}
{\it Proof.} We decompose $C$ into Whitney strata $C_i$, each
stratum fibered over $S_i\subset X$ (in $C^\infty$-topology).
Note that $\pi^{-1}Z$ is transverse to $C_i$ if and only if $Z$ is
transverse to $S_i$. (One has to stratify the set $Z$ into
Whitney strata and check transversality of each pair of strata.
For the notion of transversality of stratified cycles and their
intersections, see \cite{Go}.)
Since $X$ is homogeneous, we can move $Z$ to make it
transverse to each $S_i$. Therefore, we can assume that
$\pi^{-1}Z$ is transverse to $C$. Then we have
\begin{equation}
[C]\cdot [Z]=[C]\cdot [\pi^{-1}Z]\cdot [X]
=[C\cap \pi^{-1}Z]\cdot[X]\,.
\end{equation}
The last number is equal to the degree of the cone
class \ $z(C\cap \pi^{-1}Z, E)$ \ which is nonnegative
by \cite[Theorem 1 (A)]{FL0}.
\qed

\s

Our next aim, in the proof of Theorem \ref{main}, is to determine
the normal bundle of $LG(V)$ in
$\L(V)$. To this end, we need a general result about actions of
linear transformations on jets. Let $\cJ$ be the space of jets (of
order $k$ fixed in the beginning) of functions \ $f: (\C^m,0) \to
(\C^n,0)$, satisfying the condition \ $Df(0)=0$.

Let \ $A: \C^n\to \C^m$ be a linear map. The map $A$ acts on
$\C^m\oplus \C^n$ by
$$
(q,p)\mapsto (q+Ap,p)
$$
and it acts on $\cJ$ by transforming the graphs of the functions.
(Note that the action of $A$ is well defined because $Df(0)=0$.)
More precisely, the image of the function $f$ is the function $f_A$
satisfying the following implicit equation:
\begin{equation}\label{impl}
f(q+Af_A(q))=f_A(q)\,.
\end{equation}
We state

\begin{proposition}\label{tang} The derivative at $0$ of $A$ acting on
$\cJ$
is equal to the identity.
\end{proposition}
(We shall use this proposition in the proof of Lemma \ref{normal}, where
the setting stems from Section \ref{lag}~; in particular, we have there:
$m=n$ and $q$ (resp. $p$) are the coordinates of $W$ (resp. of $W^*$)
used in that section.)

\noindent
{\it Proof.} We compute the derivative of $A$ in the direction of an
arbitrary
$f\in\cJ$. We set
\begin{equation}
g_t=A(tf)\,.
\end{equation}
By virtue of Eq.~(\ref{impl}), the function $g_t$ is given by the implicit
equation
\begin{equation}\label{q}
tf(q+Ag_t(q))=g_t(q)\,.
\end{equation}
We want to show that \ $\frac{g_t}{t} \to f$ \ as $t\to 0$.
Since we assume
$$
g_t(0)=f(0)=0\,,
$$
it is enough to show that \ $\frac{Dg_t}{t} \to Df$ \ as $t\to 0$.
Differentiating Eq.~(\ref{q}) gives
\begin{equation}
tDf\circ(Id+A\circ Dg_t)=Dg_t\,.
\end{equation}
Hence, we get
\begin{equation}\label{dgt}
Dg_t=(Id-tDf\circ A)^{-1}\circ tDf\,.
\end{equation}
It follows from Eq.~(\ref{dgt}) that \ $\frac{Dg_t}{t} \to Df$ \ as $t\to
0$,
and the proposition has been proved.
\qed

\s

Let us come back to the setting of Section \ref{lag}, and put
$$
\L=\L(V) \ \ \ \ \hbox{and} \ \ \ \ G=LG(V)\,.
$$

Recalling that $R$ denotes the tautological vector bundle on $G$, we give
the following description of the normal bundle of $G$ in $\L$.

\begin{lemma}\label{normal}
We have a natural isomorphism
\begin{equation}
N_G\L \cong \bigoplus_{i=3}^{k+1}{{\Sym}^i}(R^*)\,.
\end{equation}
\end{lemma}
{\it Proof.} \ Let $L\in G$ be a Lagrangian linear subspace.
Let us choose a splitting
\begin{equation}\label{spl}
V\cong L\oplus L^*
\end{equation}
Using this splitting, by Lemma \ref{fib} we construct an isomorphism
\begin{equation}\label{SL}
\pi^{-1}(L) \cong \bigoplus_{i=3}^{k+1}{{\Sym}^i}(L^*)
\end{equation}
Two splittings of the exact sequence
$$
0\to L \to V \to L^*\to 0 \,,
$$
differ by a linear map $A: L^*\to L$. By Proposition \ref{tang},
the action of $A$ on the tangent space to jets is the identity.
Hence the isomorphism (\ref{SL}) does not depend of the choice
of splitting (\ref{spl}), and we have (globally)
a {\it natural} isomorphism
\begin{equation}
N_G\L\cong \bigoplus_{i=3}^{k+1}{{\Sym}^i}(R^*)\,.
\end{equation}
The lemma has been proved.
\qed

\medskip

We are now ready to complete the proof of Theorem \ref{main}.

\smallskip

Suppose that $\Sigma$ is a Lagrange singularity class (in fact, it can be
{\it any} algebraic pure-dimensional cycle in $\L$).
Let \ $i: G\hookrightarrow\L$ be the inclusion, and denote by
$$
i^*: H^*(\L,{\Z}) \to H^*(G, {\Z})
$$
the induced map on cohomology rings. We have to examine the coefficients
$\alpha_I$ of the expression
\begin{equation}
i^*[\Sigma]=\sum \alpha_I \ \Q_I(R^{*})\,.
\end{equation}
Let us fix now a strict partition $I$ with $i_1\le n$.
Invoking Proposition \ref{dual}, the coefficient $\alpha_I$ is equal to
\begin{equation}\label{coef}
i^*[\Sigma]\cdot \Omega_{I'}
\end{equation}
(intersection in $G$). Since the cohomology ring and the Chow ring
of $G$ are equal, we can compute (\ref{coef}) using the Chow groups.
Let
$$
C=C_{G\cap\Sigma}\Sigma \subset N_G\L
$$
be the {\it normal cone} of $G\cap\Sigma$ in $\Sigma$. Denote by
  \ $j: G\hookrightarrow N_G\L$ \ the zero-section inclusion.
By deformation to the normal cone \cite[\S6.1, \S6.2]{F}, we have
in $A_*G$ the equality
\begin{equation}
i^*[\Sigma]=j^*[C]\,,
\end{equation}
where $i^*$ and $j^*$ are the pull-back maps of the corresponding Chow
groups (see \cite{F}). It follows that
$$
\alpha_I=[C]\cdot \Omega_{I'}
$$
(intersection in $N_G\L$).
The bundle $R^*$ is globally generated; therefore, by Lemma \ref{normal},
the vector bundle $N_G\L$ is globally generated. The Lagrangian
Grassmannian $G=LG(V)$ is a homogeneous space with respect to the action
of the symplectic group $Sp(V)$. By Lemma \ref{int}, applied to the bundle
$N_G\L \to G$, the intersection \ $[C]\cdot \Omega_{I'}$ \ is nonnegative.

\s

This ends the proof of Theorem \ref{main}.

\s

\begin{remark}\label{rel} \rm
We may compare jets of Lagrangian manifolds enlarging the number of variables.
For a Lagrangian submanifold $L\subset V=W\oplus W^*$, we consider
$$
L\oplus \C \subset V\oplus (\C\oplus \C^*)\,.
$$
This induces an embedding $\L(W\oplus W^*)$ into $\L(W'\oplus W'^*)$,
where $W'=W\oplus \C$. In this way, we obtain  a chain of inclusions:
$$
\L(V)=\L_0 \subset\L_{1} \subset \cdots \subset \L_{r}\subset \cdots
$$
where $\L_{r}=\L((W\oplus \C^r) \oplus (W\oplus \C^r)^*)$.
We say that the singularity type $\Sigma$ is {\it closed with respect to
suspension} if one can find a sequence of singularity classes
$\Sigma_r\subset \L_{r}$ for $r\geq 0$, such that
\begin{equation}\label{stabilization}
\Sigma_r=\Sigma_{r+1}\cap\L_{r}\,
\end{equation}
and $\Sigma_0=\Sigma$.
Moreover, we assume that the property (\ref{stabilization}) holds
at the cohomological level, i.e.~the restriction
$H^*(\L_{r+1}, \Z)\to H^*(\L_{r}, \Z)$ maps the cohomology class
$[\Sigma_{r+1}]$ to $[\Sigma_r]$. (This holds when
$\Sigma_{r+1}$ and $\L_r$ intersect transversely along some smooth
open and dense subset of $\Sigma_r$.)

\s

Suppose now that $E$ is a symplectic manifold, and $E\to B$ is a
fibration with Lagrangian fibers. (Such objects have been recently
widely investigated. We do not intend here to survey this
activity, but refer the reader to, e.g., \cite{Hw} and the
references therein.) Let, in addition, $L\subset E$ be a
Lagrangian submanifold. We denote by $\rho$ the restriction of the
projection to $L$. We study the singular points of the projection
$\rho$. The definition of the singularity type (invariance with
respect to change of coordinates) allows us to define the singular
points of $\rho$ of type $\Sigma$. Denote the set of these
singular points by $\Sigma(\rho)$. Assume that the singularity
type is closed with respect to suspension (see (\ref{stabilization})).
Then, if the map $\rho$ satisfies suitable transversality conditions,
the class $[\Sigma(\rho)]\in H^*(L, {\Z})$ is equal to \ $\cT^\Sigma$
applied to the virtual bundle \ $T^*L-\rho^*T^*B$. See also
\cite[Theorem 2 and \S3.1]{Ka2}.
\end{remark}

\begin{remark}\rm

In the case of Legendre singularities, using some refinement
of our methods, we can get the following positivity result (not as transparent
as the one in the Lagrangian case). We adopt the definition of the Thom polynomial
of a Legendre singularity class from \cite{Ka2}\footnote{Our reference for Legendre
singularities and characteristic classes is \cite[\S3.2]{Ka2}.}, p.730.
This definition makes use of the {\it classifing space of Legendre singularities}
({\it loc.cit.}).
Let $\xi$ be the canonical line bundle on $BU(1)$, appearing in this definition.
Let ${\bf 1}^n$ denote the rank $n$ trivial bundle. The Thom polynomial
of a Legendre singularity class can be presented uniquely as a $\Z$-linear combination

\begin{equation}\label{alphaIj}
\sum_{j\ge 0}\sum_I \alpha_{Ij} \ \Q_I((T^*L-{\bf 1}^n)\otimes {\xi}^{\frac 1 2})\cdot t^j\,,
\end{equation}
where
   \vskip 3pt
   \hfil  $ t:={\frac 1 2}c_1(\xi^*)\in H^2(BU(1),\Z[{\frac 1 2}])\,,
    $
   \vskip 3pt
\noindent the second sum is over strict partitions $I$ and $\alpha_{Ij}\in \Z$.
We then have
$$
\alpha_{Ij} \ge 0
$$
for any strict partition $I$ and $j\ge 0$. In other words, we get
a positivity result after ``perturbing'' the argument ``$T^*L$''
by subtracting a bundle and twisting by some line bundle. If we
formally assume that $\xi$ is trivial, then we get the Lagrangian
case.
\end{remark}

\section{Examples}\label{examples}

We list now the Thom polynomials of Lagrange and Legendre singularities up
to codimension $6$. These polynomials were computed in \cite{Ka2} in
the basis of square free monomials in the Chern classes of the involved
bundles.
We present them as $\Z$-linear combinations of the
products of the ``twisted'' $\Q$-functions from (\ref{alphaIj}):
$$
\Q_I:=\Q_I((T^*L-{\bf 1}^n)\otimes \xi^{\frac 1 2})
$$
(note the difference in notation with respect to Eq.~(\ref{QI})),
and powers of $t={\frac 1 2}c_1(\xi^*)$.
The terms marked with bold give the Thom polynomials of the
corresponding Lagrange singularities.

{\parindent 0pt

\ss

   $A_2$: \ $\bf \Q_1$

   \ss

   $A_3$: \ ${\bf 3\Q_2}+ t\Q_1$

\ss

$A_4$: \ ${\bf 3\Q_{21}+12\Q_3}+ t\,10\Q_2+ t^2 2\Q_1$

\ss

$D_4$: \ $\bf \Q_{21}$

\ss

$A_5$: \ ${\bf 27\Q_{31}+60\Q_4}+t(22\Q_{21}+86\Q_3)+t^2 40\Q_2+
t^3 6\Q_1$

\ss

$D_5$: \ ${\bf 6\Q_{31}}+ t\,4\Q_{21}$

\ss

$A_6$: \  ${\bf 87 \Q_{32}+ 228 \Q_{41}+360 \Q_{5}}+
                          t(343 \Q_{31}+ 756\Q_{4})+
                          t^2(151 \Q_{21}+ 584\Q_{3})+$

$\phantom{A_6:}$ \ $          t^3 196\Q_{2}+
                                t^4 24\Q_{1}$

\ss

$D_6$: \ ${\bf 12 \Q_{32}+  24\Q_{41}}+
                          t\, 32\Q_{31}+
                           t^2 12\Q_{21}$

\ss

$E_6$: \ ${\bf 9 \Q_{32}+ 6\Q_{41}}+
                            t\,9\Q_{31}+
                             t^2 3\Q_{21}$

\ss

$A_7$: \ ${\bf 135 \Q_{321}+ 1275 \Q_{42}+ 2004\Q_{51}+  2520
\Q_{6}}+$

  $\phantom{A_6:}$ \ $t(7092\Q_{5}+ 4439 \Q_{41}+ 1713
\Q_{32})+t^2(3545\Q_{31}+ 7868\Q_{4})+ $

  $\phantom{A_6:}$ \ $t^3(1106
\Q_{21}+4292\Q_{3})+t^4 1148\Q_{2}+t^5 120\Q_{1}$

\ss

$D_7$: \
     ${\bf 24\Q_{321}+   120  \Q_{42}+ 144\Q_{51}}+
         t(152\Q_{32}+  288 \Q_{41})+
               t^2 208\Q_{31}+
                t^3 56\Q_{21}$

\ss

$E_7$: \
     ${\bf 9 \Q_{321}+ 60 \Q_{42}+  24 \Q_{51}}+
          t(56 \Q_{41}+66 \Q_{32})+
            t^2 42\Q_{31}+
                t^3 10 \Q_{21}$

\ss

$P_8$: \ $\bf \Q_{321}$.}

\bigskip

\noindent
{\bf Acknowledgements} \
We thank \"Ozer \"Ozt\"urk for pointing out some corrections in an
earlier version of the paper. We are also grateful to the referee for
a careful lecture of the manuscript and suggestion of several improvements.


\begin{thebibliography}{99}\small

\bibitem{Ar1} V. I. Arnold,
\emph{On a characteristic class entering into conditions of
quantization,}
Funk. Anal. Pril. {\bf 1} (1967), 1--14.

\bibitem{AGV} V. I. Arnold, S. M. Gusein-Zade, A. N. Varchenko,
\emph{Singularities of Differentiable Maps,}
Vol.~I, Monographs in Mathematics {\bf 82} (1985), Birkh\"auser.

\bibitem{Au2} M. Audin,
\emph{Classes caract\'eristiques lagrangiennes,}
in ``Algebraic topology, Barcelona 1986'', Springer Lect. Notes in Math.
{\bf 1298},
1--16.

\bibitem{Bo} A. Borel,
\emph{Sur la cohomologie des espaces fibres principaux et des espaces
homog\`enes de groupes de Lie compacts,}
Ann. of Math. {\bf 57} (1953), 115--207.

\bibitem{Fk} D. B. Fuks,
\emph{The Maslov-Arnold characteristic classes,}
Dokl. Akad. Nauk SSSR {\bf 178} (1968), 303--306.

\bibitem{F}  W. Fulton,
\emph{Intersection theory,} Springer, 1984.

\bibitem{FL0} W. Fulton, R. Lazarsfeld,
\emph{Positivity and excess intersections,} in ``Enumerative and
classical geometry'', Nice 1981, Progress in Math. {\bf 24},
Birkha\"user (1982), 97--105.

\bibitem{FL1} W. Fulton, R. Lazarsfeld,
\emph{The numerical positivity of ample vector bundles,}
Ann. of. Math. {\bf 118} (1983), 35--60.

\bibitem{Go} M. R. Goresky,
\emph{Whitney stratified chains and cochains,}
Trans. Amer. Math. Soc. {\bf 267} (1981), no. 1, 175--196.


\bibitem{Hw} J.-M. Hwang,
\emph{ Base manifolds for fibrations of projective irreducible symplectic
manifolds,}
arXiv:0711.3224.

\bibitem{Ka2} M. Kazarian,
\emph{Thom polynomials for Lagrange, Legendre, and critical point
function singularities,} Proc. London Math. Soc. (3) {\bf 86} (2003),
707--734.

\bibitem{L} R. Lazarsfeld,
\emph{Positivity in algebraic geometry,} Springer, 2004.

\bibitem{LP} A. Lascoux, P. Pragacz,
\emph{Operator calculus for $\Q$-polynomials and Schubert
polynomials,}
Adv. Math. {\bf 140} (1998), 1--43.

\bibitem{P1} P. Pragacz,
\emph{Algebro-geometric applications of Schur S- and Q-polynomials,}
in ``Topics in Invariant Theory'', S\'eminaire d'Alg\`ebre
Dubreil-Malliavin 1989-1990 (M.-P. Malliavin ed.),
Springer Lect. Notes in Math. {\bf 1478} (1991), 130--191.

\bibitem{PR} P. Pragacz, J. Ratajski,
\emph{Formulas for Lagrangian and orthogonal degeneracy loci;
$\Q$-polynomial approach,}
Compositio Math. {\bf 107} (1997), 11--87.

\bibitem{PW1} P. Pragacz, A. Weber,
\emph{Positivity of Schur function expansions of Thom polynomials,}
Fund. Math. {\bf 195} (2007), 85--95.

\bibitem{Va1} V. A. Vassiliev,
\emph{Characteristic classes of Lagrange and Legendre mani\-folds
that are dual to singularities of caustics and wave fronts,}
Funk. Anal. Pril. {\bf 15} (1981), 10--22.


\end{thebibliography}
\end{document}